\input amssym.def
\input amssym.tex
\documentstyle{amltd}
\begin{document}
\annalsline{157}{2003}
\received{March 21, 2000}
\startingpage{521}
\def\bye{\end{document}}
 \font\tenrm=cmr10

\def\joinrel{\mathrel{\mkern-4mu}}
\def\relbar{\mathrel{\smash-}}
\def\lrar{\relbar\joinrel\relbar\joinrel\relbar\joinrel\rightarrow}
 \def\llar{\leftarrow\joinrel\relbar\joinrel\relbar\joinrel\relbar}
\catcode`\@=11
\font\twelvemsb=msbm10 scaled 1100
\font\tenmsb=msbm10
\font\ninemsb=msbm10 scaled 800
\newfam\msbfam
\textfont\msbfam=\twelvemsb  \scriptfont\msbfam=\ninemsb
  \scriptscriptfont\msbfam=\ninemsb
\def\msb@{\hexnumber@\msbfam}
\def\Bbb{\relax\ifmmode\let\next\Bbb@\else
 \def\next{\errmessage{Use \string\Bbb\space only in math
mode}}\fi\next}
\def\Bbb@#1{{\Bbb@@{#1}}}
\def\Bbb@@#1{\fam\msbfam#1}
\catcode`\@=12

 \catcode`\@=11
\font\twelveeuf=eufm10 scaled 1100
\font\teneuf=eufm10
\font\nineeuf=eufm7 scaled 1100
\newfam\euffam
\textfont\euffam=\twelveeuf  \scriptfont\euffam=\teneuf
  \scriptscriptfont\euffam=\nineeuf
\def\euf@{\hexnumber@\euffam}
\def\frak{\relax\ifmmode\let\next\frak@\else
 \def\next{\errmessage{Use \string\frak\space only in math
mode}}\fi\next}
\def\frak@#1{{\frak@@{#1}}}
\def\frak@@#1{\fam\euffam#1}
\catcode`\@=12
\newcommand{\CA}{$C^*$-algebra}
\newcommand{\SCA}{$C^*$-subalgebra}
\newcommand{\aue}{approximate unitary equivalence}
\newcommand{\ayue}{approximately unitarily equivalent}
\newcommand{\mops}{mutually orthogonal projections}
\newcommand{\hm}{homomorphism}
\newcommand{\pisca}{purely infinite simple \CA}
\newcommand{\andeqn}{\,\,\,\,\,\, {\rm and} \,\,\,\,\,\,}
 \newcommand{\morp}{contractive completely
positive linear
 map}
\newcommand{\asmorp}{asymptotic morphism}
\newcommand{\arrow}{\rightarrow}
\newcommand{\tdsum}{\widetilde{\oplus}}
\newcommand{\pa}{\|}  
\newcommand{\ep}{\varepsilon}
\newcommand{\id}{{\rm id}}
\newcommand{\aueeps}[1]{\stackrel{#1}{\sim}}
\newcommand{\aeps}[1]{\stackrel{#1}{\approx}}
\newcommand{\dt}{\delta}
\newcommand{\yu}{\fang}
\newcommand{\ca}{{\cal C}_1}
\newcommand{\Ad}{{\rm ad}}
\newcommand{\Ik}{ {\cal I}^{(k)}}
\newcommand{\Iz}{{\cal I}^{(0)}}
\newcommand{\Ii}{{\cal I}^{(1)}}
\newcommand{\Ip}{{\cal I}^{(2)}}
\font\titr=cmr10  scaled 1200
\font\emi= cmmi10 scaled 1700 
\title{Classification of simple {\emi C}\hskip1pt\raise4pt\hbox{\titr*}-algebras and\\
higher dimensional noncommutative tori}  

 \def\titleheadline#1{\def\one{#1}\ifx\one\empty\else
\gdef\thetitle{{\frenchspacing%
\let\\ \relax
{#1}}}\fi}
\newif\ifshort
\def\shortname#1{\global\shorttrue\xdef
\theauthors{{\eightsc\uppercase{#1}}}}
\let\shorttitle\titleheadline
\shorttitle{ \eightsc\uppercase{Classification of simple} {\eightpoint \it C}
\hskip-3pt\raise2pt\hbox{{*}}\eightsc\uppercase{-algebras} }

 \acknowledgements{Research partially supported by NSF grant DMS-009790,
Zhi-Jiang Professorship
in East China Normal University and Shanghai Priority Academic Disciplines.}
 \author{Huaxin Lin}
 
\institutions{East China Normal University, Shanghai, China {\eightpoint \rm and}\\
{\eightpoint {\it Current address\/}}:  University of Oregon,
Eugene, Oregon\\
{\eightpoint {\it E-mail addresses\/}: hxlin@math.ecnu.edu.cn}\\
 \hglue.97in {\eightpoint hxlin@noether.uoregon.edu}} 
 
\vfil
\centerline{\bf Abstract}
\vglue8pt
We show that unital simple \CA s with tracial topological rank zero
which are locally approximated by  subhomogeneous \CA s can
be classified by their ordered $K$-theory.
We apply this classification result to show that certain
simple crossed products are isomorphic if they have the same
ordered $K$-theory.
In particular,
irrational higher dimensional noncommutative tori
of the form $C({\Bbb T}^k)\times_{\theta}{\Bbb Z}$ are in fact
inductive limits of circle algebras.
\vfil
\intro

  In recent years there has been rapid progress in
classification of nuclear simple \CA s. In the case that \CA s are
of real rank zero and finite, Elliott and Gong (\cite{EG}) have
proved that simple inductive limits of finite direct sums of
homogeneous \CA s (AH for brevity) of slow dimension growth with
real rank zero can be completely classified up to isomorphism by
their scaled ordered $K$-theory (with the reduction of dimension
growth proved by \cite{G2} and \cite{D}). In their remarkable paper
(\cite{EG}), they also showed that the class of  ${\rm AH}$-algebras that
they classified exhausts all possible invariants. So any general
classification theorem for simple \CA s of real rank zero and
stable rank one with weakly unperforated $K_0$ will not expand
their class. However, many interesting simple \CA s, which are
important in applications, do not arise as inductive limits of
finite direct sums of homogeneous \CA s. Therefore, it is
extremely important to have a classification theorem which covers
\CA s that are {\it not} assumed to be  ${\rm AH}$-algebras. The main
purpose of this paper is to establish such a theorem. Our general
classification result covers at least some of the well-known
interesting simple \CA s that are not known to \pagebreak  be  ${\rm AH}$-algebras.

For example, certain simple \CA s arising as dynamical systems
with minimal diffeomorphisms can be classified by their ordered
$K$-theory. More specifically, let $X_i$ ($i=1,2$) be a smooth
manifold ($X_1$ and $X_2$ may be the same) and $\sigma_i: X_i\to
X_i$ be a minimal diffeomorphism. Suppose that $(X_i, \sigma_i)$
($i=1,2$) is  uniquely ergodic. Then the resulting crossed products
$C(X_i)\times_{\sigma_i} {\Bbb Z}$ are simple \CA s with unique
traces and are isomorphic if they have the same scaled ordered
$K$-theory (which is determined by $\sigma_i$ only). A consequence
of this is that the noncommutative tori of the form $C({\Bbb T}^k)\times_{\theta} {\Bbb Z},$ where $\theta$ is an irrational
rotation, are isomorphic to  unital simple inductive limits of
circle algebras. This also generalizes an important result of
Elliott and Evans that every irrational rotation algebra is an
inductive limit of circle algebras (\cite{EE}).

\pagegoal=50pc

To classify a class of \CA s, one often needs to establish
a so-called uniqueness theorem and an existence theorem.
Uniqueness  is used to describe
two maps from one \CA\, to another as approximately equivalent
in an appropriate  sense if they carry the same $K$-theory (or $KK$-theory)
data.
Existence is often involved in
showing
that given
$K$-theory data (or $KK$-data) $\alpha,$ there is a map $\phi$ from one \CA\,
to another which carries $\alpha.$

In \cite{Ln2}, \cite{Ln3} and \cite{Ln4}, we show that a uniqueness theorem
holds for all nuclear \CA s with a reasonable and mild restriction
which, together with a ``half" existence theorem (see also
\cite{DE}), gives a number
of classification results which
do not require that  \CA s considered to be AH-algebras.
The `` half" existence theorem we mentioned above
does give us a map which carries most of
the required $KK$-data but not all.
The missing part is
the
order information from the required $KK$-data.

Suppose that $A$ is a unital simple \CA\, with $ {\rm TR}(A)=0$
(see \ref{IIIDTAF}).
Then, in \cite{Ln3}, we show  that
$A$ is always an inductive limit of $A_n,$ where each $A_n$ is
a residually finite-dimensional \CA. If $\rho_{A_n}: K_0(A_n)\to
{\rm Aff}(T(A_n))$ is the \hm\, given by
the traces of $A_n,$  and every finitely generated
subgroup of $\rho_{A_n}$ can be order-embedded into
${\Bbb Z}^k$ for some integer $k,$ then the needed existence theorem holds.
This certainly appears to be a very technical condition.
However every AH-algebra obviously
satisfies this  condition.
We show that if $A$ is an inductive limit of subhomogeneous
\CA s, then this condition is also satisfied.
In fact, a much broader class of \CA s satisfies this condition.
Combining this with our other recent results, we are able
to prove that simple \CA s with $ {\rm TR}(A)=0$ can be classified
up to isomorphism by their scaled ordered $K$-theory provided
that they are locally approximated by subhomogeneous \SCA s.
The recent results of Q. Lin and N. C. Phillips show that
every simple \CA\, arising from a dynamical system of  minimal
diffeomorphisms is in fact an inductive limit of
subhomogeneous \CA s. Therefore, the general classification
result mentioned above can be applied to those simple
\CA s when they have real rank zero.
\eject

The paper is organized as follows.
In Section 1, we revisit  ${\rm AH}$-algebras.
A few refinements of known results will be presented.
These refinements are needed for construction of some maps that
possess certain required
$KK$-data.
In Section 2, we present results which
enable us to extend some positive \hm s from ordered subgroups
of ${\Bbb Z}^k$ to some other ordered groups which are not assumed
to be divisible.
In Section 3, we use results from Sections 1 and 2 together
with the ``half" existence theorem in \cite{Ln4} (see also
\cite{DE}) to prove a full existence theorem.
We then combine our recent results to establish the main
classification theorem (\ref{IIIT1}).
A few examples of applications are presented at the end of
the paper.

\pagegoal=48pc

\demo{Acknowledgements}
We would like to acknowledge that we have benefited from
conversations with N. C. Phillips and
Guihua Gong.
\enddemo

\section{{\rm AH}-algebras revisited}

The main purpose
of this section is to prove Lemma \ref{ID}.

\proclaim{Lemma}\label{GRB}
Let $G=\lim_{n\to\infty} (G_n, \alpha_n)$ be a countable
unperforated simple ordered group with the Riesz interpolation
property{\rm ,} where $G_n=\mathbold{\oplus}_{i=1}^{l(n)}D_{i_{\phantom{1}}}^{(n)},$ $
D_i^{(n)}={\Bbb Z}$ and $\alpha_n$ is an order \hm. Suppose
that $\alpha_{m,\infty}(D_i^{(m)})\not=\{0\}$ for every $i$ and
$m.$

Then for any $0<m_1<m_2,$ there is an integer
 $N=N(m_1, m_2)>m_2$ satisfying the following\/{\rm :}
if $\pi_k\circ \alpha_{m_2, N}(G_{m_2})\not=\{0\},$
then $\pi_k\circ \alpha_{m_1, N}(D_j^{m_1})\not=\{0\}$ for $1\le j\le l(m_1)$
and
for all $k\le l(N),$ where $\pi_k: G_N\to D_k^{N}$
is the standard projection.
\endproclaim

\demo{Proof}
Let $g_i\in D_i^{(m_2)}$ be the positive generator of $D_i^{(m_2)}$
($\cong {\Bbb Z}),$ $i=1,2,\ldots ,l(m_2).$
Let $x_i=\alpha_{m_2,\infty}(g_i).$
Let $s_j$ be the positive generators for $D_j^{(m_1)}$
and let $f_j=\alpha_{m_1,\infty}(s_j),$ $j=1,2,\ldots ,l(m_1).$
From the assumption,
$f_j\not=0.$
Since $G$ is simple,
there are positive integers $k_{ij}$ such that
$k_{ij}f_j\ge x_i$ for $i=1,2,\ldots ,l(m_2).$
Therefore, there is $N=N(m_1, m_2)$ such that
$$
k_{ij}\alpha_{m_1,N}(s_j)\ge \alpha_{m_2, N}(g_i),
$$
$i=1,2,\ldots ,l(m_2).$
Fix $k\le l(N),$ define $\pi_k: G_N\to D_k^{(N)}$ to be the standard projection.
If for some $k,$ $\pi_k\circ \alpha_{m_2, N}(g_i)>0,$
then $k_{ij}\pi_k\circ \alpha_{m_1,N}(s_i)>0.$
Since $\{g_1, g_2,\ldots ,g_{l(m_2)}\}$ is a set of generators,
we see that the conclusion holds.
\enddemo

Let $G$ be a countable unperforated simple ordered group and
$T$ be the state space of $G.$
Let $\rho_G: G\to {\rm Aff}T$ be the map defined by
evaluation, i.e., $\rho_G(g)(t)=t(g)$ ($g\in G$ and $t\in T$).
It is known that
$$
G_+=\{g: \rho(g)(t)>0\,\,\,{\rm for\,\,\,all}\,\,\, t\in T\}\cup\{0\}.
$$
 
In the following lemma, it is known that one can require that $\alpha_n^{i,j}$
have multiplicity at least 2 or $\alpha_n^{i,j}=0$
(\cite{Ell2}). So the only thing   new
is that we can always assume for every $j,$
$\alpha_n^{i,j}$ has
multiplicity at least 2 (not zero) for some~$i.$

\proclaim{Lemma}\label{GROUPD}
Let $G$ be a countable unperforated simple ordered group with the
Riesz interpolation property. Suppose that ${\rm ker}\rho=0.$ Then
there are $\{G_n\},$ where $G_n$ is a finite sum of ${\Bbb Z}$
with the usual order{\rm ,} and positive \hm s $\alpha_n: G_n\to
G_{n+1}$ such that $ G=\lim_{n\to\infty} (G_n, \alpha_n).$
Furthermore{\rm ,}  nonzero $\alpha_n^{i,j}$ has
multiplicity at least $2,$ and for each $j,$ there is at least one
$\alpha_n^{i,j}\not=0,$ and for each $i,$ there is at least one
$\alpha_n^{i,j}\not=0,$ where $\alpha_n^{ i,j}: D_i^{(n)}\to
D_j^{(n+1)}$ is the partial map of $\alpha$ and $D_i^{(m)}$ is the
$i^{\rm th}$ summand of $G_m$ ($D_i^{(m)}\cong {\Bbb Z}$).
\endproclaim

\demo{Proof}
The first part of the lemma follows from \cite{EHS}.
It is the last part which needs a proof.
We write $G_n=\mathbold{\oplus}_{i=1}^{m(n)} D_i^{(n)},$ where
$D_i^{(n)}\cong {\Bbb Z}.$
Without loss of generality, we may assume that
$$
\alpha_{m,\infty}(D_i^{(m)})\not=\{0\}
$$
for all $i$ and $m.$
By \cite{Ell2}, we may also assume that $\alpha_n^{i,j}$ is either zero
or has multiplicity at least 2.

Let $G_1'=G_1$ and let
$G_2'=\mathbold{\oplus} \{D_j^{(2)}: \alpha_1^{i,j}\not=0\,\,\,{\rm for\,\,\, some }\,\,\,i\}.$
If $G_n'$ is defined, define
$$
G_{n+1}'=\mathbold{\oplus} \{ D_j^{(n+1)}: \alpha_n^{i,j}\not=0\,\,\,
{\rm for \,\,\, some}\,\,\, i\,\,\,{\rm such \,\,\,that}\,\,\, D_i^{(n)}\in G_n'\}.
$$
Set $\beta_n=(\alpha_n)|_{G_n'}.$
It is important to note that each nonzero partial map $(\alpha_n^{j,i})$
of $\beta_n$ has multiplicity at least 2 and for each $j,$ at least
one $\alpha_n^{i,j}\not=0.$

Let $G'=\lim_{n\to\infty} (G_n',\beta_n).$
It suffices to show that $G'$ is order isomorphic to $G.$
We will use $\alpha_{n,m}$ ($m>n$) for
$\alpha_m\circ \alpha_{m-1}\circ\cdots \circ \alpha_n$
and $\beta_{n,m}$ for\break $\beta_m\circ\beta_{m-1}\circ\cdots \circ \beta_n.$
Clearly $\beta_{n,m}=(\alpha_{n,m})|_{G_n'}.$

Define $\phi_1=\beta_1: G_1'\to G_2.$
By Lemma \ref{GRB}, there exists an integer $N(1, 2)>0$ such that
the conclusion of \ref{GRB} holds for $m_1=1$ and $m_2=2.$
Let $k(2)=N(1,2).$ Define $\psi_2=\alpha_{2, k(2)}.$
The conclusion of \ref{GRB} shows that $\psi_2$ maps $G_2$ to $G_{k(2)}'$
and $\psi_2\circ \phi_1=\beta_{1, k(2)}.$

Define $\phi_2=\beta_{k(2)}: G_{k(2)}'\to G_{k(2)}.$ By definition of $\psi_2,$
$\alpha_{2, k(2)}=\phi_2\circ \psi_2.$
By applying \ref{GRB} again, we obtain an integer $k(3)=N(1, k(2))$ such that
the conclusion of \ref{GRB} holds for $1$ and $k(2).$
Define $\psi_3=\alpha_{k(2), k(3)}.$
The conclusion of \ref{GRB} shows that $\psi_3$ maps $G_{k(2)}$ to
$G'_{k(3)}.$ Also, $\beta_{k(2), k(3)}=\psi_3\circ \phi_2.$
Thus we obtain the following commutative diagram:
$$
\begin{array}{clccr}
G_1'&{\stackrel{\beta_{1,k(1)}}{\longrightarrow}}&
G_{k(2)}'& {\stackrel{\beta_{k(2),k(3)}}{\longrightarrow}}& G_{k(3)}'\\[4pt]
\downarrow_{\phi_1} &  \nearrow_{\psi_2}& \downarrow_{\phi_2}&
\nearrow_{\psi_3} \\[4pt]
G_2 & {\stackrel{\alpha_{2,k(2)}}{\longrightarrow}}& G_{k(2)}&
{\stackrel{\alpha_{k(2),k(3)}}{\longrightarrow}}& G_{k(3)}.
\end{array}
$$
Continuing this construction, we obtain the next commutative
diagram:
$$
\begin{array}{clcccccr}
G_1'&{\stackrel{\beta_{1,k(1)}}{\longrightarrow}}&
G_{k(2)}'& {\stackrel{\beta_{k(2),k(3)}}{\longrightarrow}}& G_{k(3)}'&
{\stackrel{\beta_{k(3),k(4)}}{\longrightarrow}}& \cdots  {\longrightarrow} & G'\\[4pt]
\downarrow_{\phi_1} &  \nearrow_{\psi_2}& \downarrow_{\phi_2}&
\nearrow_{\psi_3} &\downarrow_{\phi_3} &  \nearrow_{\psi_4} &\\[4pt]
G_2 & {\stackrel{\alpha_{2,k(2)}}{\longrightarrow}}& G_{k(2)}&
{\stackrel{\alpha_{k(2),k(3)}}{\longrightarrow}}& G_{k(3)} &
{\stackrel{\alpha_{k(3),k(4)}}{\longrightarrow}} &\cdots  {\longrightarrow} &G.
\end{array}
$$

Therefore $G\cong G'.$ Since each $\phi_n$ and $\psi_n$ is positive,
this isomorphism is in fact an order isomorphism.
\enddemo

\numbereddemo{Definition}\label{SPA}
 Let $f: S^1\to S^1$ be a degree $k$ map ($k>1$), i.e., a continuous map
with winding number $k.$ We let (following 4.2 in \cite{EG})
$T_{II,k}=D^2\mathbold{\cup}_f S^1,$ the finite connected CW complex obtained
by attaching a $2$-cell $D^2$ to $S^1$ via the map $f.$
Note that  $K_0(C(T_{II,k}))={\Bbb Z}\oplus {\Bbb Z}/k{\Bbb Z}$ and
$K_1(C(T_{II,k}))=\{0\}.$
Let $g: S^2\to S^2$ be a degree $k$ map ($k>1$).
Let $T_{III,k}=D^3\mathbold{\cup}_g S^2$ be the connected finite CW complex
obtained by attaching a $3$-cell $D^3$ to $S^2$ via the map $g.$
Note that $K_0(C(T_{III,k}))={\Bbb Z}$ and $K_1(C(T_{III,k}))={\Bbb Z}/k{\Bbb Z}$
(see 4.2 in \cite{EG}). 
\enddemo

\numbereddemo{Definition}\label{IPE}
{\rm Let $C=PM_n(C(X))P,$ where $P\subset M_n(C(X))$
is a projection with rank $r(x)$ at point $x.$
Note that if $X$ is connected, $r(x)$ is a constant.
Let $B$ be  another \CA.
A map $\omega: C\to B$ is said to be a {\it point-evaluation},
if $\omega=h\circ \pi_{x},$ where $x\in X$ is a point,
$\pi_{x}(f)=f(x)$ maps $C$ to $M_{r(x)}$ and $h: M_{r(x)}\to B$
is a \hm. Suppose that $e\in M_{r(x)}$ is a minimal
projection. We say that $\omega$
is a point-evaluation with a minimal projection $h(e).$
}
\enddemo

The following is a refinement of a result in \cite{EG}.
Only part (3) is new.

\proclaim{Theorem}\label{ICUT}
For any countable simple weakly unperforated scaled ordered
group $(G, G_+, [u])$ with the Riesz interpolation property and
any countable abelian group $F,$ there exists a
unital simple \CA\, $A$ of real rank zero with the following properties\/{\rm :}
\begin{itemize}
\item[{\rm (1)}] $A=\lim_{n\to\infty}(A_n, h_n),$
where each $A_n$ is a finite direct sum of\break
$P_iM_{m(i)}(C(X_i))P_i,$ where
$X_i=Y_1\vee Y_2\vee \cdots Y_m$ and
each $Y_i=S^1, S^2, T_{II, k},\break T_{III, l},$
or a point\/{\rm ;}

\item[{\rm (2)}] $(K_0(A), K_0(A)_+, [1_A])=(G, G_+, u)$ and $K_1(A)=F$ and

\item[{\rm (3)}] ${\rm ker}\rho_{K_0(A)}=
\lim_{n\to\infty} ({\rm ker}\rho_{K_0(A_n)}, (h_n)_{*}).$
\end{itemize}

\endproclaim

\demo{Proof} Let $G_0=\rho_G(G).$ Then $G_0$ is an unperforated
ordered group.
By Lemma \ref{GROUPD}, we may write
$G_0=\lim_{n\to\infty} (G_n,\alpha_n),$
where $D_m=\mathbold{\oplus}_{i=1}^{l(m)} {\Bbb Z}$ with the usual order and
$\alpha_n^{(i,j)}: D_i^{(n)}\to D_j^{(n+1)}$ \pagebreak has order at least $2$
and satisfies the rest of the requirements of \ref{GROUPD}.
Let ${\rm ker}\rho_G=\lim_{n\to\infty}(H_n,\beta_n)$
and let $F=\lim_{n\to\infty}(F_n, \gamma_n),$ where
$H_n$ and $F_n$ are  finitely generated abelian groups.
Let ${\tilde G_n}=\alpha_{n,\infty}(G_n).$ We may write
${\tilde G_n}\subset {\tilde G_{n+1}},$ $n=1,2,\ldots\, .$

Let ${\tilde S_n}$ be the subgroup of $G$ generated by ${\rm ker}\rho_G$ and
$G_n'$
such that $\rho_G(G_n')={\tilde G_n}.$
Since ${\tilde G}_n$ is free, we may write
${\tilde S_n}={\tilde G}_n\oplus {\rm ker}\rho_G.$
Let $\imath_n: {\tilde S_n}\to {\tilde S}_{n+1}$ be the embedding.
Then $\imath_n(g\oplus h)=g\oplus (\imath_n^{(0)}(g)+h),$
where $\imath_n^{(0)}: {\tilde G}_n\to {\rm ker}\rho$
is a \hm\, (note that if $G=G_0\oplus {\rm ker}\rho,$ then
one may choose $\imath_n^{(0)}=0$).
Set $S_n=G_n\oplus {\rm ker}\rho_G.$
Define $j_n: S_n\to S_{n+1}$ by
$j_n(g\oplus h)=\alpha_n(g)\oplus (\alpha_n^{(0)}(g)+h),$
where $\alpha_n^{(0)}=\imath_n^{(0)}\circ \alpha_{n,\infty}.$
Then the following is commutative:
$$
\begin{array}{ccc}
{\tilde S_n}&{\stackrel{\imath_n}{\to}}&{\tilde S}_{n+1}\\
\uparrow_{\alpha_{n,\infty}} & & \uparrow_{\alpha_{n+1,\infty}}\\
S_n & {\stackrel{j_n}{\to}}& S_{n+1}.\\
\end{array}
$$
Set $\eta_n=(j_n)|_{G_n}$ and $\dt_n=\eta_n\oplus \beta_n.$
Set $$(G_n\oplus H_n)_+=\{(g,x): g>0, x\in H_n\}\cup\{0\}.$$
Note
that $(G, G_+)=\lim_{n\to\infty}(G_n\oplus H_n, (G_n\oplus H_n)_+,\dt_n).$

Let $X_{n}=Y_1\vee Y_2\vee \cdots \vee Y_{t(n)},$
where $Y_i=S^1, S^2, T_{II, k}$ or $T_{III, l},$ such
that $K_0(C(X_{n}))={\Bbb Z}\oplus H_n$ with
\begin{eqnarray*}
 {\rm ker}\rho_{C(X_n)}&=&H_n,\\
 K_0(C(X_{n}))_+&\subset &\{(z,x):  z>0, x\in H_n\}\cup \{0\},\\
 \{(y,x): y\ge 3,x\in H_n\}&\subset& K_0(C(X_{n}))_+
\end{eqnarray*}
(see 4.17 in \cite{EG}) and
$K_1(C(X_{n}))=F_n.$
Suppose that $u=\dt_n(u_n)$ and that $\pi_j(u_n)\ge 4,$
where $\pi_j: G_n\to {\Bbb Z}$ is the projection to the $j^{\rm th}$ coordinate.

Let $A_1= P_1M_{s(1)}(C(X_{1}))P_1\oplus B_1,$ where
$B_1=\mathbold{\oplus}^{l(1)-1}M_2,$ $P_1\in\break M_{r(1)}(C(X_{1,1}))$ is a
projection so that $[P_1]=\pi_1(u_1)$ ($s(1)>\pi_1(u_1)$). We have
\begin{eqnarray*}
K_0(A_1)&=&G_1\oplus H_1,\\
{\rm ker}\rho_{A_1}&=&H_1,\\
 K_0(A_1)_+&\subset&\{(g,x): g\in (G_1)_+\}\oplus \{0\},\\
 \{(g,x):g\in 3(D_1^{(1)})_+\}&\subset& K_0(A_1)_+,\\
\noalign{\noindent and}
 K_1(A_1)&=&F_n.
\end{eqnarray*} Set $C_1=P_1M_{r(1)}(C(X_{1}))P_1.$  Denote by
$r(1,1)$ the rank of $P_1$ ($r(1,1)\le r(1)$), and let $r(1,i)=2$
for $2\le i\le l(1).$

Let $N$ be the required integer in Lemma 3.27 in \cite{EG}
corresponding to the space $X_{1}$ and $\ep<1/4.$
Using \ref{GROUPD},
choose $k_2$ such that each nonzero partial map $\alpha_{1,k_2}^{i,j}$
has multiplicity $N+3.$ Since, for each $j,$
$\alpha_{1,k_2}^{i,j}\not=0$ for some $i,$
$\pi_j(u_{k_2})\ge N+3$ for all $j.$
Let $P_2$ be a projection in $M_{r(2)}(C(X_{k_2}))$
for some integer $r(2)$ so that
$[P_2]=\pi_1(u_{k(2)})$ and $Q_2\le P_2$ be a projection
with rank $=$ the multiplicity of $\alpha_{1,k_2}^{1,1}$
times $r(1,1).$
Since the rank of $Q_2$ is at least $(N+3)r(1,1),$ we may assume that
$P_2-Q_2$ is a trivial projection.
Set $B_2 =\mathbold{\oplus}_{i=1}^{l(k_2)-1} M_{r(2,i)},$
where $r(2,i)=\sum_{j=1}^{l(1)} m_{1,j,i}r(1,i)$ and $m_{1,j,i}$
is the multiplicity of $\alpha_{1,k_2}^{j,i}.$
Set $C_2=P_2M_{r(2)}(C(X_{k_2}))P_2.$

Define $h_{1,2}:C_1\to B_2$ by
$$
h_{1,2}(f)=\sum_{j=1}^{l(k_2)-1}\omega_{1,j}(f),
$$
where $\omega_{1,j}: C_1\to M_{r(2,j)}$ is the  point evaluation
at a point $x_1\in X_1$ with  minimal projection
having rank $m_{1,j,i}$ (see \ref{IPE}).
Define
$h_{1,1}: B_1\to B_2$ according to the multiplicity
($(l(1)-1)\times (l(k_2)-1)$) matrix
$(\alpha_{1,k_2}^{i,j})_{i=2,j=1}^{l(1)-1,l(k_2)-1}.$
Define
$h_{2,1}: B_1\to C_2$ according to the multiplicity
of $\alpha_{1,k_2}^{i,1}$ so that the minimal projections
are trivial.

We define $h_{2,2}: C_1\to C_2$ by
applying 3.27 in \cite{EG}.
Write $Q_2=Q_{2,0}\oplus Q_{2,1}$ with
$Q_{2,1}$ having rank $12\times r(1,1)$ and
$Q_{2,0}$ a trivial projection.
Define $h_{2,2}: C_1\to C_2$ by
$h_{2,2}=\phi_{1,2}\oplus \phi_{0,2},$ where
$$
\phi_{0,2}(f)=\left( \begin{array}{ccccc} \omega_{x_1}(f) & & & & \cr
                       &      \omega_{x_2} (f) &&&\cr
                       & & \ddots && \cr
                       &&& \omega_{x_{K_1}}(f) & \cr
                       &&&& \omega_{x_1}'(f) \end{array}\right),
$$
where $\{x_1,\ldots ,x_{K_1}\}$ is $1/4$-dense in $X_1,$
$\omega_{x_i}$ is the point-evaluation at $x_i$
such that the minimal projection
is a trivial projection in $C_2$ of rank $12,$ $i=1,2,\ldots ,K_1,$
and  $\omega_{x_1}'$ is the point-evaluation at $x_i$
such that  the minimal projection is
a trivial projection and $\omega_{x_1}(1_{C_1})
=Q_{2,0}-\sum_{i=1}^{K_1}\omega_{x_i}(1_{C_1}).$
Furthermore,
$\phi_{1,2}: C_1\to Q_{2,1}C_2Q_{2,1}$ is given
by the map $\phi_1$ as described in the proof
of 3.27 in \cite{EG} so that
$\phi_{1,2}(1_{C_1})=Q_{2,1},$
$[\phi_{1,2}]|_{H_1}=\beta_{1,k_2}$ and $[\phi_{1,2}]_{K_1(C_1)}=
\gamma_{1,k_2}$ ($K_1(C_1)=F_1$).

Define $\Phi_{1,2}: A_1\to A_2$ by
$$
\Phi_{1,2}=\left( \begin{array}{cc} h_{1,1} & h_{1,2}\cr
                  h_{2,1} & h_{2,2}\end{array}\right)
$$
according to the  decomposition $A_1=B_1\oplus C_1$ and $A_2=B_2\oplus C_2.$
From the above construction,
we verify that
$$
[\Phi_{1,2}]|_{K_0(A_1)}=\dt_{1,k_2}\andeqn [\Phi_{1,2}]|_{K_1(A_1)}=\gamma_{1,k_2}.
$$
By continuing this construction,
we obtain $\{A_n\}$ and $\Phi_{n, n+1}.$
Then $A=\lim_{n\to\infty} (A_n, \Phi_{n,n+1})$ has
real rank zero
and $(K_0(A),K_0(A)_+, [1_A])=(G, G_+, u)$
(see the proof of 4.18 in \cite{EG} for example)
and
$K_1(A)=F.$ The simplicity of~$G$ also implies
that $A$ is simple. This also follows the fact
that the multiplicity of each $\alpha_n^{(i,j)}$ is
at least $4.$
So (1) and (2) follow. Part (3) also follows
from the fact that $[\Phi_{1,2}]|_{K_0(A_1)}=\dt_{1,k_2}.$
\enddemo

\numbereddemo{Definition}\label{KKL}
{\rm Let $C_n$ be a commutative \CA\, with $K_0(C_n)={\Bbb Z}/n{\Bbb Z}$
and $K_1(C_n)=0.$ Suppose that $A$ is a \CA.
Then $K_i(A, {\Bbb Z}/k{\Bbb Z})=K_i(A\otimes C_k).$
Let ${\bf P}(A)$ be the set of all projections in
$M_{\infty}(A),$ $M_{\infty}(C(S^1)\otimes A),$
$M_{\infty}((A\otimes C_m{\tilde)})$ and
$M_{\infty}((C(S^1)\otimes A\otimes C_m{\tilde )}).$
We have the following commutative diagram (\cite{Sc}):
$$
\begin{array}{ccccc}
K_0(A) &\to &K_0(A, {\Bbb Z}/k{\Bbb Z})
&\to& K_1(A)\\
\uparrow_{\bf k} & & & &\downarrow_{\bf k}\\
K_0(A) & \leftarrow & K_1(A, {\Bbb Z}/k{\Bbb Z})
& \leftarrow & K_1(A).
\end{array}
$$
As in \cite{DL}, we use the notation
$$
{\underline K}(A)= \bigoplus_{i=0,1, n\in {\Bbb Z}_+}
K_i(A;{\Bbb Z}/n{\Bbb Z}).
$$
By
${\rm Hom}_{\Lambda}({\underline K}(A),{\underline K}(B))$
we mean all \hm s from ${\underline K}(A)$ to ${\underline K}(B)$
which respect the direct sum
decomposition and the so-called Bockstein operations (see \cite{DL}).
It follows from \cite{DL} that if $A$ satisfies the Universal
Coefficient Theorem, then ${\rm Hom}_{\Lambda}({\underline K}(A),{\underline K}(B))
=KL(A,B).$


Let $A$ and $B$ be two \CA s and $L: A\to B$ a completely positive
linear map. Then $L$ induces maps from $A\otimes C_m$ to $B\otimes C_m,$
from $C(S^1)\otimes A\otimes C_m$ to $C(S^1)\otimes B\otimes C_m,$
namely, $L\otimes {\rm id}.$ For convenience, we will also
denote the induced map by $L$.
Let $A$ and $B$ be \CA s, let $L: A\to B$ be a \morp, let $\ep>0$ and let
${\cal F}\subset A$ be a subset.  Now, $L$ is said to be
${\cal F}$-$\ep$-multiplicative,
if
$$
\|L(xy)-L(x)L(y)\|<\ep
$$
for all $x,y\in {\cal F}.$ Given a projection $p\in {\bf P}(A),$
if $L$ is ${\cal G}$-$\ep$-multiplicative with sufficiently large
${\cal G}$ and sufficiently small $\ep,$ $L(p)$ is close to a
projection. Let $L(p)'$ be that projection. Fix  a finite subset
${\cal P}_1\subset {\bf P}(A).$ It is easy to see that  $L(p)'$
and $L(q)'$ are in the same equivalence class of projections of
${\bf P}(A),$ if $p$ and $q$ are in ${\cal P}_1$ and are in the
same equivalence class of projections of ${\bf P}(A),$ provided
that ${\cal F}$ is sufficiently large and $\ep$ is sufficiently
small. We use $[L](p)$ for the class of projections containing
$[L](p)'.$ In what follows, whenever we write $[L](p),$ we assume
that ${\cal F}$ is sufficiently large and $\ep$ is sufficiently
small so that $[L](p)$  is  well-defined on ${\cal P}_1.$
Furthermore, abusing the language, we write $[L]([p])$ as well as
$[L](p),$ where $[p]$ is the equivalence class containing $p.$

Suppose that $q$ is in ${\bf P}(A)$ with
$[q]=k[p]$ for some integer $k$; by adding sufficiently many
elements (partial isometries) in ${\cal F},$ we can
assume that $[L](q)=k[L](p).$
Suppose that $G$ is a finitely generated group generated
by ${\cal P}$ and $G={\Bbb Z}^n\oplus {\Bbb Z}/k_1{\Bbb Z}\oplus
\cdots {\Bbb Z}/k_m{\Bbb Z}.$
Let $g_1,g_2,\ldots ,g_n$ be free generators of ${\Bbb Z}^n$
and $t_i\in {\Bbb Z}/k_i{\Bbb Z}$ be the generator with order $k_i,$
$i=1,2,\ldots ,m.$ Since every element in $K_0(C)$ (for any unital \CA\, $C$) may be
written as $[p_1]-[p_2]$ for projections $p_1, p_2\in A\otimes M_l,$
for some $l>0,$
with sufficiently large ${\cal F}$ and
sufficiently small $\ep,$ one can define
$[L](g_j)$ and $[L](t_i).$
Moreover (with sufficiently large ${\cal F}$ and
sufficiently small $\ep$), the order of $[L](t_i)$ divides $k_i.$
Then we can define a map $[L]|_{G}$
by
defining $[L](\sum_i^n n_ig_i+\sum_j^mm_jt_j)=
\sum_i^k n_i[L](g_i)+\sum_j^mm_j[L](t_j).$
Thus $[L]$ is a group \hm\, on $G.$
Note, in general, that  $[L]|_{{\cal P}}$ may {\it not}
coincide with $[L]|_{G}$ on ${\cal P}.$
However, if ${\cal F}$ is large enough and $\ep$ is
small enough, they coincide.
In what follows,
if ${\cal P}$ is given, we say $[L]|_{G}$ is well-defined  and write $[L]|_{G}$
if $[L]|_{\cal P}$ is well-defined, $[L]|_{G}$ is well-defined
and is a \hm\, and $[L]|_{{\cal P}}=[L]|_{G}$ on ${\cal P}.$
}
\enddemo

\numbereddemo{Definition}\label{IBB}
{\rm We denote by ${\cal C}$ the family of
all unital simple \CA s of real rank zero
which are direct limits of finite direct sums
of unital hereditary \SCA s of $M_n(C(X))$ (for various $n$),
where $X$ is a connected finite CW complex of dimension
no more than $3$ (and $X$ may be different in the sums).
This is precisely the class
of simple \CA s classified in \cite{EG}.
}
\enddemo

\proclaim{Lemma}\label{ID}
Let $A$ be a unital simple \CA \, in ${\cal C}.$
Let $G_0$ be a finitely generated   subgroup
of $K_0(A)$ with decomposition
$G_0=G_{00}\oplus G_{01},$ where $G_{00}\subset {\rm ker}\rho_A$ and
$G_{01}$ is a finitely generated free group such that
$(\rho_A)|_{G_{01}}$ is injective. Suppose that
${\cal P}\subset {\underline{K}}(A)$ is a finite subset
which generates a subgroup $G$ such that
$G\cap K_0(A)=G_0.$

Then{\rm ,} for any $\ep>0,$ any finite subset
${\cal F}\subset A,$
any $1>r>0,$ and any integer $K,$
there is an ${\cal F}$\/{\rm -}\/$\ep$\/{\rm -}\/multiplicative
map $L: A\to A$ satisfying the following\/{\rm :}\/
\begin{itemize}
\item[{\rm (1)}] $[L]|_{\cal P}$ and $[L]|_{G}$ are  well\/{\rm -}\/defined and
$[L]|_{G}$ is positive on $G,$

\item[{\rm (2)}] $[L]|_{G\cap {\rm ker}\rho_A}={\rm id}|_{G\cap {\rm
ker}\rho_A},$ $[L]|_{G\cap K_0(A, {\Bbb Z}/k{\Bbb Z})}={\rm
id}| _{G\cap K_0(A, {\Bbb Z}/k{\Bbb Z})},$
$[L]|_{G\cap K_1(A)}\break ={\rm id}|_{G\cap K_1(A)}$ and 
$[L]|_{G\cap K_1(A, {\Bbb Z}/k{\Bbb Z})}=
{\rm id}|_{G\cap K_1(A, {\Bbb Z}/k{\Bbb Z})}$
for those $k$ with
$
G\cap K_i(A, {\Bbb Z}/k{\Bbb Z})\not=\emptyset
$ {\rm (}$i=0,1${\rm ),}

\item[{\rm (3)}] $|\rho_A\circ [L](g)|\le r|\rho_A(g)|$ for all $g\in G\cap
K_0(A).$

\item[{\rm (4)}] Let $g_1, g_2,..,g_l$ be positive generators of
$G_{01}.$ Then{\rm ,} there are $f_1,\ldots ,f_l\in K_0(A)_+$ such that
$$
g_i-[L](g_i)=Kf_i, \,\,\,i=1,2,\ldots ,l.
$$
\end{itemize}
\endproclaim

\demo{Proof}
We may write
$A=\lim_{n\to\infty}(A_n, \phi_{n, n+1}),$
where $A_n=\mathbold{\oplus}_{i=1}^{m(n)}B_{n,i},$
each $$B_{n,i}=P_{n,i}M_{J(n,i)}(C(X_{n,i}))P_{n,i}$$
for some connected finite CW complex $X_{n,i}$ with ${\rm dim}X_{n,i}\le 3$
and $P_{n,i}\in M_{J(n,i)}(C(X_{n,i}))$ is a projection, and
$\phi_{n,n+1}: A_n\to A_{n+1}$ is a \hm, as constructed
in \ref{ICUT}.
In particular, (3) in \ref{ICUT} holds.

Let $\phi_{n,i,j}: B_{n,i}\to B_{n+1,j}$ be the partial map
determined by $\phi_{n, n+1}.$
Fix $0<r<1.$
By the proof of  \ref{ICUT}, we may assume that
$$
\phi_{n,i,j}(f)={\rm diag}({\tilde \phi_{n,i,j}}(f), \psi_{n,i,j}(f)),
$$
where ${\tilde \phi_{n,i,j}}(1_{B_{n,i}})$ has rank
no more than $12\times {\rm rank}(1_{B_{n,i}})$ and
$$
\psi_{n,i,j}(f)=\sum_{s=1}^{l(n,i,j)}\omega_{n,i,j,s}(f),
$$
where $\omega_{n,i,j,s}$
is a point-evaluation (see \ref{IPE})
such that $\omega_{n,i,j,s}(1_{B_{n,i}})=q_{n,i,j,s}$, 
$\{q_{n,i,j,s}\}$ is a set of mutually orthogonal projections
in $B_{n+1,j}$ and $\{q_{n,i,j,s}\}$ are
equivalent (trivial) projections in $B_{n+1, j}$ with rank
at least $12\times {\rm rank}(1_{B_{n,i}})$
for $s=2,\ldots ,l(n,i,j),$ and $q_{n,i,j,1}$ has rank
${\rm rank 1_{B_{n+1,j}}}-12\times l(n,i,j){\rm 1_{B_{n,i}}}.$

By choosing larger $n,$ we may assume that
${\cal F}\subset \phi_{n,\infty}(A_n).$
Furthermore,
we may also assume that
$$
G\subset [\phi_{n,\infty}](\underline{K}(A_n)).
$$
Let $G'$ be a finitely generated subgroup of $\underline{K}(A_n)$ such that
$[\phi_{n,\infty}](G')\break=  G.$
To save notation without loss of generality, we may assume
that $G_0=[\phi_{n,\infty}](K_0(A_n)).$
Write $K_0(A_n)=F_0\oplus F_1,$
where $F_0={\rm ker}\rho_{A_n}.$
By (3) of \ref{ICUT}, we may assume, without loss of
generality,  that
$[\phi_{n,\infty}](F_0)=G_{00}.$
Let $K_0$ be the integer such that
$$
G'\cup K_i(A, {\Bbb Z}/k{\Bbb Z})=\{0 \}
$$
whenever $k\ge K_0$ and $i=0,1.$
Put $K_1=K(K_0)!.$
By replacing $n+1$ by a larger integer, if necessary, we may assume that
$$
l(n,i,j)>(K_1+1)K_1(1/r)\,\,\,\,\,{\rm for\,\,\, all}\,\,\,
n,i,j.
$$

Let
$$
l(n,i,j)-1=K_1v(n,i,j)'+r(n,i,j),
$$
where $v(n,i,j)$ and $r(n,i,j)$ are nonnegative integers with
$r(n,i,j)<K_1.$

Define $$\Phi_{n,i,j}(f)={\rm diag}({\tilde \phi_{n,i,j}}(f), \psi_{n,i,j}'(f)),\hbox{ where }
\psi_{n,i,j}'(f)=\sum_{s=1}^{1+r(n,i,j)}\omega_{n,i,j,s}(f)
$$
and define
$$
\Phi_{n,i,j}'(f)=\sum_{s=2+r(n,i,j)}^{l(n,i,j)}\omega_{n,i,j,s}(f)
\hbox{ 
for all $f\in B_{n,i}.$}
$$
Define $\Phi(f)=\mathbold{\oplus}_{i,j}\Phi_{n,i,j}(f)$
and $\Phi'(f)=\mathbold{\oplus}_{i,j}\Phi_{n,i,j}'(f)$ for $f\in A_n.$
Since $\Phi'$ has finite-dimensional range,
we know that
$$
[\Phi']|_{K_0(A_n)\cap {\rm ker}\rho_{A_n}}=0
\andeqn
[\Phi']|_{K_1(A_n)}=0.
$$
For each $i$ and $j,$ $\Phi_{n,i,j}'(f)$ is a direct sum
of $K_1v(n,i,j)$ many point-evaluations.
Since $X_{n,i}$ is connected and
$(K_0)!|K_1,$ by considering each $[\Phi_{n,i,j}'],$
we conclude that (for $k=2,\ldots ,K_0$)
$$
[\Phi']|_{K_0(A_n, {\Bbb Z}/k{\Bbb Z})}=0
\andeqn
[\Phi']|_{K_1(A_n, {\Bbb Z}/k{\Bbb Z})}=0.
$$
 
Therefore
$$
\begin{array}{rlrl}
 [\Phi]|_{F_0}&\hskip-6pt =[\phi_{n,n+1}]|_{F_0}, & 
[\Phi]|_{K_0(A_n, {\Bbb Z}/k{\Bbb Z})}
&\hskip-6pt =[\phi_{n,n+1}]|_{K_0(A_n, {\Bbb Z}/k{\Bbb Z})},\\[5pt]
 \left[\Phi\right]|_{K_1(A_n)}&\hskip-6pt =[\phi_{n,n+1}]|_{K_1(A_n)} &  
\andeqn
\left[\Phi\right]|_{K_1(A_n, {\Bbb Z}/k{\Bbb Z})}
&\hskip-6pt =\left[\phi_{n,n+1}\right]|_{K_1(A_n, {\Bbb Z}/k{\Bbb Z})}
\end{array}
$$
for every $k=1,2,\ldots , K_0.$ 
Thus
\begin{eqnarray*}
 [\phi_{n+1,\infty}\circ \Phi]|_{F_0}&\hskip-6pt=\hskip-6pt&[\phi_{n,\infty}]|_{F_0},\\[5pt]
[\phi_{n+1,\infty}\circ \Phi]|_{K_0(A_n, {\Bbb Z}/k{\Bbb Z})}
&\hskip-6pt=\hskip-6pt&[\phi_{n,\infty}]|_{K_0(A_n, {\Bbb Z}/k{\Bbb Z})},\\[5pt]
\left[\phi_{n+1,\infty}\circ \Phi\right]|_{K_1(A_n)}&\hskip-6pt=\hskip-6pt&[\phi_{n,\infty}]|_{K_1(A_n)}\\
\noalign{\noindent and}
[\phi_{n+1,\infty}\circ \Phi]|_{K_1(A_n, {\Bbb Z}/k{\Bbb Z})}
&\hskip-6pt=\hskip-6pt&[\phi_{n,\infty}]|_{K_1(A_n, {\Bbb Z}/k{\Bbb Z})}
\end{eqnarray*}
for every $k=1,2,\ldots , K_0.$
It is clear that 
$$
{\rm ker}\phi_{n+1,\infty}\circ \Phi=
{\rm ker}\phi_{n,n+1}.
$$
(In fact, we could assume that
${\rm ker}\phi_{n+1,\infty}\circ \Phi={\rm ker}\phi_{n,n+1}=\{0\}.$)
Thus, $\phi_{n+1,\infty}\circ \Phi$ induces
a map $\Psi: \phi_{n,\infty}(A_n)\to A.$
Set $A_n'=\phi_{n,\infty}(A_n).$
Denote by ${\imath}: A_n'\to A$ the
embedding. From above,
we have
$$
\begin{array}{rlrl}
 [\Psi]|_{{\rm ker}\rho_{A_n'}}&\hskip-6pt =[\imath]_{{\rm ker}\rho_{A_n'}},&
[\Psi]|_{K_1(A_n')}&\hskip-6pt=[\imath]_{K_1(A_n')},\\[5pt]
 [\Psi]|_{K_0(A_n', {\Bbb Z}/k{\Bbb Z})}
&\hskip-6pt=[\imath]_{K_0(A_n', {\Bbb Z}/k{\Bbb Z})},&
\andeqn
[\Psi]|_{K_1(A_n', {\Bbb Z}/k{\Bbb Z})}
&\hskip-6pt=[\imath]_{K_1(A_n', {\Bbb Z}/k{\Bbb Z})}
\end{array}
$$
for $k=1,2,\ldots\, .$

Since $A$ is unclear, for any finite subset ${\cal F}_1\subset A$
and $\dt>0,$ there exists (see 4.1 in \cite{Ln4}) a completely
positive linear map $L: A\to A$ such that
$$
\|L(a)-\Psi(a)\|<\dt
$$
for all $a\in {\cal F}_1.$
We may assume that ${\cal F}\subset {\cal F}_1$ and
$\dt<\ep.$
Furthermore, by choosing even larger ${\cal F}_1$ and
small $\dt,$
we may assume that
$$\begin{array}{rlrl}
 [L]|_{G_{00}}&\hskip-6pt={\rm id}|_{G_{00}},&
[L]|_{G\cap K_1(A)}&\hskip-6pt={\rm id}|_{G\cap K_1(A)}\\[5pt]
 [L]|_{G\cap K_0(A, {\Bbb Z}/k{\Bbb Z})}
&\hskip-6pt={\rm id}|_{G\cap K_0(A, {\Bbb Z}/k{\Bbb Z})}&
\andeqn
[L]|_{G\cap K_1(A, {\Bbb Z}/k{\Bbb Z})}
&\hskip-6pt={\rm id}|_{G\cap K_1(A, {\Bbb Z}/k{\Bbb Z})}
\end{array}
$$
for all $k$ so that $G\cap  K_i(A, {\Bbb Z}/k{\Bbb Z})\not=\emptyset$
($i=0,1$).
From  $l(n,i,j)\ge\break (K_1+1)K_1(1/r)$
we conclude that
$$
\tau(L(a)))<r\tau(a)
$$
for all $a\in A$ and $\tau\in T(A).$
This implies that $L$ satisfies (1), (2)  and (3).
For $s\ge 2,$ $[\omega_{n,i,j,s}]=[\omega_{n,i,j,2}].$
Since, for any $z\in K_0(B_{n,i}),$
$$
[\phi_{n,i,j}](z)-[\Phi_{n,i,j}](z)
=K(K_0)!v(n,i,j)[\omega_{n,i,j,2}](z),
$$
(4) also follows.
\enddemo

\section{Extensions of positive \hm s on ordered groups}

Let $G$ be a group, $G_0$ be a subgroup of $G$ and
$F$ be another group. In general, to extend a \hm\, $\phi: G_0\to F$
to a \hm\, ${\tilde \phi}: G\to F$ requires $F$ to be divisible.
If $G$ is an ordered group and $\phi$ is positive,
much more is required of $F$
to obtain a positive extension.
In this section, we will present a positive extension theorem
which neither requires divisibility nor  completeness of $F.$
This result is rather special but is essential for us
to construct maps in \ref{IIIL2}.

\proclaim{Lemma}\label{IIL1}
Let $G\subset {\Bbb Z}^k$ be an ordered subgroup
and let $$e_1=(1,0,\ldots ,0), \ldots ,e_k=(0,\ldots ,0,1)\in {\Bbb Z}^k.$$
Suppose that $S\subset \{1,2,\ldots ,k\}$ such that
if $i\in S,$
then there exists a positive integer $m_i$ such that $m_ie_i\in G$
and if $i\not\in S,$ then $me_i\not\in G$ for any $m\in {\Bbb Z}\setminus \{0\}.$
Then for any positive \hm\, $\phi: G\to {\Bbb R}$ with
$\phi(G_+\setminus \{0\})\subset {\Bbb R}_+\setminus \{0\},$
$$
{\rm inf}\{\phi(g)/n: g\in G, g\ge ne_i\}>\sup \{\phi(g)/m: g\in G, g\le me_i\}>0
$$
for all $i\not\in S.$
\endproclaim

\demo{Proof}
This follows from Lemma 2.10 in \cite{Ln6}.
\enddemo

\proclaim{Lemma}\label{IIL2}
Let $G\subset {\Bbb Z}^k,$
let $$e_1=(1,0,\ldots ,0),\ldots ,e_k=(0,\ldots ,0,1)\in {\Bbb Z}^k,$$
and $S\subset \{1,2,\ldots ,k\}$ be as in Lemma {\rm \ref{IIL1}.}
Let $T$ be  a Choquet simplex and  $F\subset {\rm Aff}(T)$ be
a dense subgroup such that $f\in F_+\setminus \{0\},$
$f(t)>0$ for all $t\in T.$
Suppose that $\phi: G\to F$ is a positive \hm .
Let
$$
U_j(t)=\inf\{\phi(g)(t)/n: n\in {\Bbb N}, g\in G, ne_j\le g\},
$$
$$
L_j(t)=\sup\{\phi(g)(t)/n: n\in {\Bbb N}, g\in G, g\le ne_j\}
$$
and
$H_j(t)=U_j(t)-L_j(t).$
Then
$$
\liminf_{t\to t_0}H_j(t)>0
$$
for every $t_0\in T,$ if $j\not\in S.$
\endproclaim

\demo{Proof}
Clearly $H_j(t)\ge 0$ for all $t\in T.$
Suppose that $j\not\in S.$
By Lemma~\ref{IIL1},
if $H_j(t)=0$ for some $t\in T,$ then
$j\in S.$
Therefore $H_j(t)>0$ for every $t\in T.$

We extend $\phi$ on ${\Bbb Q}G$ by defining
$\phi(rg)=r\phi(g)$ for all $r\in {\Bbb Q}$ and $g\in G.$

Note that, since ${\Bbb Q}G$ is finite-dimensional, if
$\{x_n\}\subset {\Bbb Q}G$ is a bounded sequence, so is
$\{\phi(x_n)\}.$ Note also that  $H_j(t)>0$ for all  $t\in T.$
Without loss of generality, to simply notation, we may assume that
$j=1.$  There are $g_n\in {\Bbb Q}G$ with $g_n\ge e_1,$
$g_n=(1, r_1^{(n)}, r_2^{(n)},\ldots ,r_k^{(n)}),$ where $r_i^{(n)}\in
{\Bbb Q}_+^{\phantom{|^1}}$ for $i=2,3,\ldots ,k$ and $t_n\in T$ such that $t_n\to
t_0$ and $\phi_n(g_n)(t_n)(t)\to \liminf_{t\to t_0}U_j(t);$ and,
there are $y_n\in {\Bbb Q}G$ with $y_n\le e_1,$
$y_n=(1,q_1^{(n)},\ldots ,q_k^{(n)}),$ where $-q_i^{(n)}\in {\Bbb Q}_+$ for $i=2,3,\ldots ,k$ and $s_n\in T$ such that $s_n\to t_0$ and
$\phi(y_n)(s_n)\to \limsup_{t\to t_0}L_j(t).$

By 2.9 in \cite{Ln6},
for each $t\in T,$
there exist $\alpha_i>0,$ $i=1,2,\ldots ,k$, such that
$\phi(z)(t)=\langle z,\omega\rangle$ for $z\in G,$ where $\omega=(\alpha_1,\ldots ,\alpha_k).$
Therefore $\{r_i^{(n)}\}$ is a bounded sequence for every $i=2,\ldots ,k.$
Similarly, $\{q_i^{(n)}\}$ is a bounded sequence for every $i=2,\ldots ,k.$
Thus  $\{\phi(g_n)\}$ and $\{\phi(y_n)\}$ are (uniformly)  bounded (on $T$).
Since ${\Bbb Q}G$ is finite-dimensional,
$\{\phi(g_n)\}$ and $\{\phi(y_n)\}$ are pre-compact subsets of ${\rm Aff}(T).$
Thus, without loss of generality, we may assume that
$\phi(g_n)\to g$ and $\phi(y_n)\to x$ uniformly on $T$, where $g,\, x\in {\rm Aff}(T).$
Since each $g_n\ge e_1$ and $y_n\le e_1,$  we conclude that
$g(t)\ge  U_j(t)>0$ and  $x(t)\le L_j(t)$ for all $t\in T.$
Furthermore,
$$
g(t_0)=\liminf_{t\to t_0}U_j(t)\andeqn
x(t_0)=\limsup_{t\to t_0}L_j(t).
$$

We assume that $H_j(t)>0$ for all $t\in T.$
Therefore $g(t)>x(t)$ for all $t\in T.$
Since $x$ and $g$ are continuous  and $T$ is compact,
$$
{\rm inf}\{g(t)-x(t): t\in T\}>0.
$$
This implies that
\vglue8pt
\hfill ${\displaystyle 
\liminf_{t\to t_0}H_j(t)>0
\,\,\,\,\,{\rm for\,\,\, all},\,\,\, t_0\in T.
}$
\enddemo

\proclaim{Lemma}\label{IIL3}
In Lemma {\rm \ref{IIL2},} if $K>0$ is a previously given integer{\rm ,} then
there exists $f\in F\subset {\rm Aff}(T)$ such that
$$
L_j(t)<Kf<U_j(t)
$$
for all \pagebreak $t\in T.$
\endproclaim

\demo{Proof}  
Since $T$ is compact and $H_j$ is upper semi-continuous,\break
$\liminf_{t\to t_0}H_j(t)  >0$ for all $t_0\in T$
 and the fact that $H_j(t)>0$ for
all $t\in T$  implies that
$
\inf\{H_j(t): t\in T\}>0.
$
Let
$
a=(1/32)\inf\{ H_j(t):t\in  T\}>0.
$
Then
$
0<31a< \liminf_{t'\to t} H_j(t')$  for  all 
$t\in T.$
Fix $t_0\in T,$ let $g$ and $x$ be as in the proof of \ref{IIL2}.
Then
\begin{eqnarray*}
\limsup_{t'\to t_0}L_j(t')&\hskip-6pt<\hskip-6pt&
x(t_0)+a/8 <x(t_0)+a/4\\ &\hskip-6pt<\hskip-6pt&g(t_0)-a/2<g(t_0)-a/4 
<\liminf_{t'\to t_0}U_j(t').
\end{eqnarray*}
We have
$$
x(t')\le L_j(t')\andeqn U_j(t')\le g(t')
$$
for any $t'\in T.$ So, in particular,
$$
x(t')\le g(t')
$$
for all $t'\in T.$
Therefore there is a neighborhood $O(t_0)$ such that the following holds:
\begin{eqnarray*}
L_j(t')& \hskip-6pt<\hskip-6pt& \limsup_{t''\to t}L_j(t'')+a/16<x(t')+a/8<x(t')+a/4<g(t')-a/2\\
& \hskip-6pt<\hskip-6pt&g(t')-a/4<\liminf_{t''\to t}U(e_1)(t'')-a/8<U(e_1)(t') 
\end{eqnarray*}
for all $t'\in O(t_0).$

Since $T$ is compact, there are $O(t_1), O(t_2),\ldots ,O(t_l),$
such that
$\mathbold{\cup}_{i=1}^l O(t_l)\break \supset T.$
Note that (with  $x_i$ and $g_i$ corresponding to $O(t_i)$)
$x_i(t)+a/4, g_i-a/2\in {\rm Aff}(T)$, $i=1,2,\ldots ,l.$
Set
\begin{eqnarray*}
{\check x}&=&(x_1+a/4)\vee (x_2+a/4)\vee \cdots (x_l+a/4)
\\
\noalign{\noindent and}  {\hat g}&=&(g_1-a/2)\wedge (g_2-a/2)\wedge \cdots (g_l-a/2).
\end{eqnarray*}
Since $x_i\le L_j$ and $U_j\le g_i$ for $i,j=1,2,\ldots , l,$
$
{\check x}\le {\hat g}.
$
Since $T$ is a Choquet simplex, ${\rm Aff}(T)$ has the Riesz interpolation
property (see II.3.11 in \cite{Alf}). Thus
there is $h\in {\rm Aff}(T)$ such that
$$
{\check x} \le h\le {\hat g}.
$$
Therefore (by considering each $O(t_i)$), we have
$$
L_j(t)+a/16<{\check x}(t) \le h(t)\le {\hat g}(t)
<U_j(t)-a/16
$$
for all $t\in T.$
Hence,
$$
L_j(t)<h(t)-a/32<h(t)<U_j(t)-a/16.
$$
Since $F$ is dense in ${\rm Aff}(T),$ there exists
$f\in F$ such that
$$
\|f-(1/K)h\|< a/128K.
$$
Therefore,
$$
L_j(t)<h(t)-a/32<Kf(t)<h(t)+a/32<U_j(t) 
$$
for all $t\in T.$
\enddemo

\proclaim{Lemma}\label{IIL5}
Let $G$ and $D$ be ordered groups.
Suppose that there is a surjective \hm\, $\rho: G\to D$ such
that $g\ge 0$ if and only if $\rho(g)\ge 0.$
Then{\rm ,} for any other ordered group $G'$ and a \hm\, $\psi:
G'\to G,$ $\psi$ is positive if and only if
$\rho\circ \psi$ is positive.
\endproclaim

{\it Proof}. This is evident.
\hfill\qed

\proclaim{Lemma}\label{IIL4}
Let $G\subset {\Bbb Z}^k$ be an ordered subgroup{\rm ,}
$T$ be a Choquet simplex and $F\subset {\rm Aff}(T)$ be an
ordered dense subgroup with the strict ordering
{\rm (}\/i.e.{\rm ,} $f\in F_+\setminus \{0\}$
implies $f(t)>0$ for all $t\in T).$
Suppose that $G'$ is an ordered group with a
surjective map $\rho: G'\to F$ such that
$g\in G'_+$ if and only if $\rho(g)\in F_+.$
There exists an integer $K$ depending only on $G$  satisfying
the following\/{\rm :}
Suppose that $\phi: G\to G'$ is a positive \hm\,
such that $\phi(g)>0$ for all $g\in G_+\setminus \{0\}$ and such that
$\phi(g_i)=Kf_i$ for the generating set $\{g_1,\ldots ,g_m\}\subset
G,$ where $f_i\in G';$
 then there is a  positive \hm\, ${\tilde \phi}: {\Bbb Z}^k\to G'$
such that ${\tilde \phi}|_G=\phi$ and
$\rho\circ {\tilde \phi}(e_j)\in G'_+\setminus \{0\},$ where
$e_1=(1,0,\ldots ,0),\ldots ,e_k=(0,\ldots ,0,1)$ are the standard
generators of ${\Bbb Z}^k.$
\endproclaim

\demo{Proof} Let $P_j: {\Bbb Z}^k\to {\Bbb Z}^j$ be the
projection on the first $j$ coordinates. Let $G_j$ be the subgroup
generated by $P_j({\Bbb Z}^k)$ and $G,$ $j=1,2,\ldots ,k.$ Set
$G_0=G.$ Let $S_0$ be the subset of $\{1,2,\ldots ,k\}$ such that
there is a positive integer $m_t$ with $m_te_t\in G_0$ whenever
$t\in S_0.$
 We may
assume that $S_0=\{1,2,\ldots ,l\}$ ($l\ge 0$). Let
$$
I_{ji}=\{m\in {\Bbb Z}: m_ie_i\in G_j\}.
$$
Then $I_{ji}$ is a subgroup of ${\Bbb Z}.$
Let
$m_{ji}={\rm min}\{|m|\in I_{ji}\setminus \{0\}\}$ and
$$
K_j=\prod_{i\in S_j} m_{ji}.
$$
Set $J=\prod_{j=0}^kK_j$ and $K=J^k.$

This integer does not depend on $\phi$ but only on $G.$
Let $\{g_1,g_2,\ldots ,g_n\}$ be a generating set for $G.$
Suppose that there are $f_1,\ldots ,f_n\in G$ such that
$\phi(g_i)=Kf_i,$ $i=1,2,\ldots\, .$

The condition that $\phi(g_i)\!=\!Kf_i$ for some $f_i\!\in\! G'$ implies
that $(1/K)\phi(m_ie_i)\break\in G'$ for some positive $m_i$ and for all
$i\le l.$
 So we define ${\tilde
\phi}(e_i)=(1/m_i)\phi(m_ie_i)$ ($i=1,2,\ldots ,l$). Note that
$(1/J^{k-1}){\tilde \phi}(g)\in G'_+$ for $g\in G_+.$ Since  (for
$j\le l$)
\begin{eqnarray*}
\sup\{\rho\circ\phi(g)/m: g\in G, m>0&&\andeqn g\le me_j\}>0,
\\
\sup\{\rho\circ{\tilde \phi}(g)/m: g\in G_l, m>0&&\andeqn g\le me_j\}>0.
\\
\noalign{\noindent
So $\rho\circ{\tilde \phi}$ maps
$(G_l)_+\setminus \{0\}\subset F_+\setminus \{0\}.$}
\end{eqnarray*}

\vglue-18pt
It follows from \ref{IIL3} that there is $f\in F_+$ such that
$$
L_{l+1}(t)<Kf(t)<U_{l+1}(t)
$$
for all $t\in T,$
where
$$
L_{l+1}(t)=\sup\{\rho\circ{\tilde \phi}(g)/m:m\in {\Bbb N},
g\in G_l \andeqn g\le me_{l+1}\}
$$
and
$$
U_{l+1}(t)=\inf\{\rho\circ{\tilde \phi}(g)/m: m\in {\Bbb N},
g\in G_l \andeqn g\ge me_{l+1}\}.
$$
Let $g\in G'$ such that $\rho(g)=f.$ Define ${\tilde
\phi}(e_{l+1})=Kg.$ It follows from \cite{GH} that ${\tilde \phi}$
extends $\phi$ and $\rho\circ {\tilde \phi}$ is positive. By
\ref{IIL5}, ${\tilde \phi}$ is positive. Furthermore, $\rho\circ
{\tilde \phi}$ maps $(G_{l+1})_+\setminus \{0\}$ into
$F_+\setminus \{0\}.$ Let $S_{1}$ be a subset of $\{1,2,\ldots ,k\}$
such that $m_te_t\in G_{l+1}$ for some positive integer $m_t$
whenever $t\in S_{1}.$ Note that for any $g\in (G_{l+1})_+,$
$(1/J^{k-1}){\tilde \phi}(g)\in G'_+.$ If $l+2\in S_2,$ define
${\tilde \phi}(e_{l+2})=(1/m_{l+1}){\tilde \phi}(m_{l+1}e_{l+2}).$
Otherwise, by \ref{IIL3}, there is $f_{l+2}\in F_+$ such that
$$
L_{l+2}(t)<Kf_{l+2}<U_{l+2}(t)
$$
for all $t\in T.$
Choose $g_{l+2}\in G'$ such that $\rho(g_{l+2})=f_{l+2}.$
Define ${\tilde \phi}(e_{l+2})=Kg_{l+2}.$
In either case we obtain a positive extension ${\tilde \phi}$ on
$G_{l+2}^{\phantom{|^|}}$ and $(1/J^{k-2}){\tilde \phi}(g)\break\in F_+$ for all
$g\in (G_{l+1})_+.$
So by an induction argument we obtain ${\tilde \phi}:
{\Bbb Z}^k\to G'$ as desired.
\enddemo

\vglue-9pt
\section{A classification
theorem for simple nuclear \CA s\\ with tracial topological
rank zero} 
\vglue-8pt

\numbereddemo{Definition}\label{III1DBD}
{\rm A \CA\, $A$ is said be in ${\cal BD}$ if there is an integer $k>0$
such that every irreducible representation of
$A$ is finite-dimensional and its dimension is no more than $k.$
The integer $k$ is called the bound.
A\break \CA\, $A$ is said be in ${\cal LBD}$ (locally ${\cal BD}$)
if for any $\ep>0$ and any finite subset ${\cal F}\subset A,$
there exists $B\in {\cal BD}$ such
that
$$
{\rm dist}(x, B)<\ep\,\,\,\,\,\,\,\,{\rm for \,\,\,all}\,\,\, x\in {\cal F}.
$$}
\enddemo

\proclaim{Lemma}\label{IIIL0} \hskip-9pt
Let $A$ be a unital separable \CA\, in {\rm BD}
with the bound~$k,$
let $1=f_1,f_2,\ldots ,f_m\in \rho_A(K_0(A)_+)$ and let $G$ be the subgroup generated
by $f_1,\ldots ,f_m.$
Then there exists finite\/{\rm -}\/dimensional irreducible
representations $\pi_1, \pi_2,\ldots ,\pi_N$ of $A$ such that
$$
g\mapsto ({\rm tr}\circ(\pi_1)_*(g),\ldots , {\rm tr}\circ(\pi_N)_*(g)))
$$
is an order embedding from $G$ to ${\Bbb Q}^N,$ where
${\rm tr}$ is the normalized trace on matrix algebras.
\endproclaim

\numbereddemo{{R}emark}\label{IIIR0}
In Lemma \ref{IIIL0}, set $x={\rm max}\{{\rm tr} \circ \pi_i([1_A]): i=1,\ldots ,N\}.$
Then the composition
$$
g\mapsto ({\rm tr}\circ(\pi_1)_*(g),\ldots , {\rm tr}\circ(\pi_N)_*(g))
\mapsto x   ({\rm tr}\circ(\pi_1)_*(g),\ldots , {\rm tr}\circ(\pi_N)_*(g))
$$
gives an order embedding from $G$ into ${\Bbb Z}^N.$
\enddemo

\demo{Proof of Lemma {\rm \ref{IIIL0}}}
This is proved in \cite{Ln7}.
Denote by $s(A)$ the set of normalized traces on $A$ defined
by $t(a)={\rm tr}\circ \pi(a)$ ($a\in A$), where $\pi$ is a
finite-dimensional irreducible representation of $A$ and
${\rm tr}$ is the standard normalized trace on matrix algebras, equipped
with the weak*-topology. It follows from Corollary 2.7 in
\cite{Ln7} that the closure of the convex hull of $s(A)$ is
$T(A),$ the tracial space of $A.$ Therefore the map
$f\mapsto f|_{s(A)}$ is an order embedding from ${\rm Aff}(T(A))$
to $C(s(A)).$ Let $D_A: K_0(A)\to C(s(A))$ be defined by
$p\mapsto p(t)$\break $(t\in s(A)$), where $p$ is a projection in
$A\otimes {\cal K}.$ Then $D_A$ is a positive \hm.
Therefore the map from $\rho_A(K_0(A))$ ($ \subset {\rm Aff}(T(A))$)
to $D_A(K_0(A)_+)$ defined by restriction
is an order isomorphism.
The lemma then follows from Lemma 2.4 in \cite{Ln7}.
\enddemo

\proclaim{Lemma}\label{IIIL2}
Let $A$ be a unital \CA \, in ${\cal LBD}$ which is
a simple \CA\, with stable rank one and weakly unperforated
$K_0(A)$ and let $B\in {\cal C}$ which is a unital
simple \CA\, with real rank zero.
Suppose that $\alpha\in {\rm Hom}_{\Lambda}({\underline K}(A),
\underline{K}(B))$ which
gives an order isomorphism from
$(K_0(A), K_0(A)_+${\rm ,} $[1_A], K_1(A))$
to $(K_0(B), K_0(B)_+, [1_B], K_1(B)).$
Then there is a sequence of \morp s
$L_n: A\to B$ such that{\rm ,} for any finite subset
${\cal P}\subset {\bf P}(A),$
$$
[L_n]|_{\cal P}=\alpha|_{\cal P}
$$
for all sufficiently large $n$ and
$$
\|L_n(ab)-L_n(a)L_n(b)\|\to 0
$$
for all $a,\, b\in A.$
\endproclaim

\demo{Proof}
Fix a finite subset ${\cal F}\subset A$ and a finite
subset ${\cal P}\subset {\bf P}(A).$
Since $A$ is in ${\cal LBD},$
we may assume that there exists $A_m\subset A$ such that
$A_m\in {\cal BD},$ ${\cal F}\subset A_m$ and
${\cal P}\subset [j]({\cal G}),$
where ${\cal G}\subset {\bf P}(A_m)$ is a finite subset and
$j: A_m\to A$ is the embedding.
Set $\alpha=\beta\circ [j]\in {\rm Hom}_{\Lambda}({\underline{K}}(A_m),
({\underline{K}}(B))_+.$  Since both $A$ and $B$ are simple,
and $\alpha|_{K_0(A)}$ is an order isomorphism, for any $g\in
K_0(A)\setminus \{0\},$ $\rho_B\circ \beta(g)\not=0.$
Note that $A_m$ satisfies the Universal
Coefficient Theorem.
It follows from 5.9 in \cite{Ln4} that there exist
a sequence of \morp s $\psi_n: A_m\to B\otimes {\cal K}$ and
a \hm\, $H_n: A_m\to B\otimes {\cal K}$ with finite-dimensional
range such that
\begin{equation}\label{IIIeq1}
[\psi_n]|_{\cal G}=\beta|_{\cal G}+[H_n]|_{\cal G}. \speqnu{{\rm e}3.1}
\end{equation}
Let $G'=K_0(A_m)\cap G,$ where $G=G({\cal G})$
is the finitely generated subgroup generated by ${\cal G}.$
We may assume that there are projections
$p_1,\ldots ,p_l\in M_k(A_m)$ for some $k$ such that
$G'$ is generated by $[p_1],\ldots ,[p_l].$
We may assume that (\ref{IIIeq1}) above holds for $G'.$
Let $G_0=\rho_{A_m}(G').$
It follows from \ref{IIIL0} (and \ref{IIIR0})
 that $G_0\subset K_0(\pi(A_m))\subset {\Bbb Z}^k$ for
some integer $k>0,$ where $\pi: A_m\to C$ is a
(surjective) \hm\, from $A_m$ to a finite-dimensional
\CA\, $C.$
Let $K$ be the integer associated with $G_0$ defined by \ref{IIL4}.

Let $K_1$ be the integer such that
$G\cap K_0(A, {\Bbb Z}/k{\Bbb Z})=\emptyset$
for all $k>K_1.$

Let $\Psi_n=\psi_n\oplus H_n\oplus\cdots \oplus H_n,$ the direct sum
of $K(K_1)!-1$ copies of $H_n.$
Thus
$$
[\Psi_n]|_{\cal G}=\alpha|_{\cal G}+K(K_1)![H_n]|_{\cal G}.
$$
If $F$ is a finite-dimensional \CA\, then one has the   following
commutative diagram:
$$ \begin{array}{ccccc}
  K_0(F)  &\lrar&
 K_0(F,{\Bbb Z}/k{\Bbb Z})&\lrar &K_1(F)\\
\big\uparrow&&&&\big\downarrow\\
K_0(F)&\llar&
 K_1(F, {\Bbb Z}/k{\Bbb Z})& \llar &K_1(F) $$
\end{array}
$$ 
where $K_0(F,{\Bbb Z}/k{\Bbb Z})=K_0(F)/kK_0(F),$
$K_1(F)=0,$ $ K_1(F, {\Bbb Z}/k{\Bbb Z})=0.$
Since $H_n$ factors through a finite-dimensional \SCA,
it is easy to check that
$$
[H_n]|_{K_1(A)\cap G }=0,\,\,\, [H_n]|_{K_1(A, {\Bbb Z}/k{\Bbb Z})\cap G}=0
\andeqn [H_n]|_{{\rm ker}\rho_A(K_0(A))\cap G}=0.
$$
Moreover,
$$
(K_1)![H_n]|_{K_0(A, {\Bbb Z}/k{\Bbb Z})\cap G}=0\,\,\,(k\le K_1).
$$
Therefore
\begin{eqnarray*}
 [\Psi_n]|_{K_1(A)\cap G}&=&\alpha|_{_1(A)\cap G },\\
\left[\Psi_n\right]|_{K_1(A, {\Bbb Z}/k{\Bbb Z})\cap G}&=&
\alpha|_{K_1(A, {\Bbb Z}/k{\Bbb Z})\cap G }  \\  
\left[\Psi_n\right]|_{{\rm ker}\rho_A(K_0(A))\cap G}&=&\alpha|_{{\rm ker}\rho_A(K_0(A))\cap G} \\
\noalign{\noindent 
and}
[\Psi_n]|_{K_0(A, {\Bbb Z}/k{\Bbb Z})\cap G}&=&
\alpha|_{K_0(A, {\Bbb Z}/k{\Bbb Z})\cap G}.\end{eqnarray*}
Choose $r>0$ such that $r<{1\over{JK(K_1)!+1}},$
where $J$ is an integer so that $[H_n(1_A)]\le [1_{M_J(B)}].$
Set $G_1'=[\Psi_n]({\cal G}).$
Let ${e_i}$ be projections such that
$[{\bar e_i}]=[\Psi_n]([p_i])-[{\bar q_i}]$ ($=\alpha([p_i])$), where
$[{\bar q_i}]=K(K_1)![H_n]([p_i]),$ $i=1,2,\ldots ,l.$
Set $G_1=[\Psi_n]{\cal G})\cup\{[{\bar e_i}], \alpha([p_i]), [{\bar q_i}], i=1,\ldots ,l\}.$

Let $L: B\to B$ be as described in \ref{ID} associated
with $G_1$ (with $\ep>0$ and ${\cal F}$ to be determined
later).
Set $\Phi_n=L\circ \Psi_n.$ Then
\begin{eqnarray*}
 [\Phi_n]|_{K_1(A_m)\cap G}&=&\alpha|_{_1(A_m)\cap G },\\
\left[\Phi_n\right]|_{K_1(A_m, {\Bbb Z}/k{\Bbb Z})\cap G}&=&
\alpha|_{K_1(A_m, {\Bbb Z}/k{\Bbb Z})\cap G }\\
\left[\Phi_n\right]|_{{\rm ker}\rho_{A_m}(K_0(A_m))\cap G}&=&\alpha|_{{\rm ker}\rho_{A_m}(K_0(A_m))\cap
G}\end{eqnarray*}
and
\begin{eqnarray*}
 [\Phi_n]|_{K_0(A_m, {\Bbb Z}/k{\Bbb Z})\cap G}&=&
\alpha|_{K_0(A_m, {\Bbb Z}/k{\Bbb Z})\cap G}.\end{eqnarray*}

We also have
$$
\rho_B\circ [\Psi_n](g)\le r\rho_B\circ \alpha(g)\,\,\,
{\rm  for}\,\,\, g\in K_0(A)\cap G\,\,
$$
and
$$\alpha([p_i])-[\Phi_n]([p_i])=K(K_1)![p_i],\,\,\,\,\,i=1,\ldots ,l,
$$
 where
$f_i\in K_0(B).$
Note that $(\alpha-[\Psi_n])(g)>0$ for all  $G_0\setminus \{0\},$ since
$r<1/(KJ(K_1)!+1).$ Note that also
with the order embedding $G_0\subset K_0(\pi(A))\break\subset {\Bbb Z}^k$
and the choice of $K,$ by \ref{IIL4}, there is a positive
\hm\, $\Psi:  K_0(\pi(A_m)) \to K_0(B)$ such that
$$
\Psi|_{G_0}=(\alpha-[\Phi_n])|_{G_0}.
$$
Since $B$ has real rank zero and stable rank one,
we obtain a \hm\, $h_n: \pi(A_m)\to (1-p)B(1-p)$
such that $[h_n]=\Psi,$ where $p=1_B-\Phi_n(1_A)$
(we may assume that $\Psi_n(1_A)\le 1_B,$ since $r<1/2(KJ(K_1)!+1)$).
We set $L_n=\Phi_n\oplus h_n\circ\pi$ with sufficiently small
$\ep$ (depends on $n$) and the finite subset ${\cal F}_1$ (which
is larger than $\Phi_n({\cal F})$).
Now
\vglue9pt \hfill $
[L_n]|_{\cal P}=\alpha|_{\cal P}. $
\enddemo
 
\numbereddemo{Definition}\label{IIIDTAF}
{\rm Recall (\cite{Ln6}) that a unital simple \CA\, $A$ has
tracial topological rank zero (written $ {\rm TR}(A)=0$) if, for any
$\ep>0,$  any finite subset ${\cal F}$ and any nonzero positive
element $a,$ there exists a nonzero projection $p\in A$ and a
finite-dimensional $C^*$-subalgebra $B\subset A$ with $1_{B}=p$
such that
\begin{itemize}
\item[(1)] $\|px-xp\|<\ep$ for all $x\in {\cal F},$

\item[(2)] $pxp\in_{\ep} B$ for all $x\in {\cal F}$
and

\item[(3)] $1-p$ is equivalent to $q\in {\overline{aAa}}.$
\end{itemize}

In \cite{Ln3}, a simple unital \CA\, $A$ with $ {\rm TR}(A)=0$ is called
TAF. It is shown in \cite{Ln3} that  a simple \CA\, $A$ with $ {\rm TR}(A)=0$
is quasidiagonal and has real rank zero, stable rank one and
weakly unperforated $K_0(A).$\break
Every simple AH-algebra $A$ with slow dimension growth and real rank zero
has $ {\rm TR}(A)=0$ (this was proved in \cite{EG}. See 2.6 in \cite{Ln3}).

One could prove the following result directly.
But this follows from the main result in \cite{Ln5}.}
\enddemo

\proclaimtitle{cf.\ \cite{Ln5}}
\proclaim{Theorem}\label{UNI}
Let $A\in {\cal LBD}$ be a unital separable simple
\CA\, with  unique normalized trace. Suppose that $A$ has real rank zero{\rm ,}
stable rank one and weakly unperforated $K_0(A).$
Then $ {\rm TR}(A)=0.$
\endproclaim

\numbereddemo{Definition}\label{IIIapp}
Let $L_1,\, L_2: A\to B$ be two linear maps, $\ep>0$ and
${\cal F}\subset A$ be a subset.
We write
$$
L_1\approx_{\ep} L_2 \,\,\,\,\,\,{\rm on}\,\,\,{\cal F}
$$
if $\|L_1(x)-L_2(x)\|<\ep$ for all $x\in {\cal F}.$
\enddemo

\proclaimtitle{Theorem 2.3 in \cite{Ln4}}
\proclaim{Theorem} \label{IIIL3}
Let $A$ be a separable unital nuclear simple
\CA\, with $ {\rm TR}(A)=0$ satisfying the {\rm UCT.}
Then{\rm ,} for any $\ep>0,$ and any finite subset ${\cal F}\subset A,$
there exist $\dt>0,$ a finite subset ${\cal P}\subset
P(A)$ and a finite subset ${\cal G}\subset A$ satisfying the following\/{\rm :}
for any unital \CA \, $B$ of real rank zero and stable rank one
with weakly unperforated $K_0(B),$ and
any two ${\cal G}$-$\dt$-multiplicative morphisms $L_1, L_2: A\to B$
with
$$
[L_1]|_{\cal P}=[L_2]|_{\cal P},
$$
there exists a unitary $U\in B$ such that
$$
{\rm ad}(U)\circ L_1\approx_{\ep} L_2
\,\,\,{\rm on}\,\,\,{\cal F}.
$$
\endproclaim

\proclaim{Theorem}\label{IIIT1}
Let $A$ and $B$ be two unital \CA s in ${\cal LBD}$
with $ {\rm TR}(A)= {\rm TR}(B)=0$ satisfying the {\rm UCT.}
Suppose that there is an
order isomorphism $\alpha:
(K_0(A), K_0(A)_+, [1_A], K_1(A))\to
(K_0(B), K_0(B)_+, [1_B], K_1(B)),
$
then there is an isomorphism $h: A\to B$ such
that $h_*=\alpha.$
\endproclaim 

\demo{Proof}
Since $A$ satisfies the UCT, there is an (invertible)
element $z\in KK(A,B)$ such
that $z|_{K_i(A)}=\alpha.$
We will use $\alpha$ for the corresponding element
in $KL(A,B).$

Fix a dense sequence $\{x_n\}$ of
the unit ball of $A$ and a dense sequence $\{y_n\}$
of the unit ball of $B.$
Set $\ep_n=1/2^n.$
Let ${\cal P}_1={\cal P}(\ep_1/2, \{x_1\})\subset  {\bf P}(A),$
$\dt_1=\dt(\ep_1/2, \{x_1\})>0$ and ${\cal F}_1
={\cal G}(\ep_1/2, \{x_1\})$ be
as in \ref{IIIL3}  associated with $\ep_1/2>0$ and
finite subset $\{x_1\}.$
We assume that $x_1\in {\cal F}_1.$
By \ref{IIIL2}, there is a \morp\,
$L_1: A\to B$ such that
$$
\|L_1(ab)-L_1(a)L_1(b)\|<\dt_1/2
$$
for all $a, b\in {\cal F}_1$ and
$
[L_1]|_{{\cal P}_1}
=\alpha|_{{\cal P}_1}.
$
Let ${\cal F}_1'=L_1({\cal F}_1)\cup \{y_1\}.$
Let ${\cal Q}_1={\cal P}(\ep_2/2, {\cal F}_1')\subset {\bf P}(A),$
${\cal G}_1={\cal G}(\ep_2/2,{\cal F}_1')\subset B$
and $d_1=\dt(\ep_2/2,{\cal F}_1')>0$ be
as in \ref{IIIL3} (for $B$) associated with $\ep_2/2$ and
${\cal F}_1'.$
We may assume that ${\cal Q}_1\supset [L_1]({\cal P}_1),$
${\cal G}_1\supset {\cal F}_1'$ and $d_1<{\rm min}\{\ep_2/2, \dt_1/2).$

By \ref{IIIL2}, there exists
$\Psi_1': B\to A$ such that
$$
\|\Psi_1'(cd)-\Psi_1'(c)\Psi'(d)\|<d_1/2
$$
for all $c, d\in {\cal F}_1'$ and
$[\Psi_1']|_{{\cal Q}_1}=\alpha^{-1}|_{{\cal Q}_1}.$
Then $\Psi_1'\circ L_1$ is
$\dt_1$-multiplicative on ${\cal F}_1$ and
$[\Psi'_1\circ L_1]|_{{\cal P}_1}=[\id]_{{\cal P}_1}.$
It follows from \ref{IIIL3}
that there is a unitary $u_1\in A$ such that
$$
\Ad u_1\circ (\Psi'_1\circ L_1)\approx_{\ep_1/2} \id_A\,\,\,\,\,\,{\rm on}\,\,\,
\{x_1\}.
$$
Define $\Psi_1=\Ad u_1\circ \Psi_1'.$
Set ${\cal S}_2=\Psi_1({\cal G}_1)\cup\{x_1, x_2\}.$
Let ${\cal F}_2={\cal G}(\ep_2/2, {\cal S}_2),$
${\cal P}_2={\cal P}(\ep_2/2, {\cal S}_2)$ and
$\dt_2=\dt(\ep_2/2, {\cal S}_2)>0$ (for $A$) associated
with $\ep_2/2$ and ${\cal S}_2.$
We may assume that
${\cal F}_2\supset {\cal S}_2,$ ${\cal P}_2\supset
[\Psi_1]({\cal Q}_1)$ and $\dt_2<{\rm min}\{\ep_2/2,d_1/2).$

It follows from \ref{IIIL2} that there exists
$L_2': A\to B$ such that
$$
\|L_2'(ab)-L_2'(a)L_2'(b)\|<\dt_2/2
$$
for all $a, b\in {\cal F}_2$ and
$$
[L_2']|_{{\cal P}_2}=\alpha|_{{\cal P}_2}.
$$
Note that $\alpha(\alpha^{-1})=[\id].$
By \ref{IIIL3}
there is a unitary $v_2\in B$ such that
\vglue4pt
\centerline{${\displaystyle
\Ad v_1\circ L_2'\circ \Psi_1\approx_{\ep_2/2} \id_B
\,\,\,\,\,\,{\rm on}\,\,\,{\cal G}_1.
}$}
\vglue2pt\noindent 
Set $L_2=\Ad v_1\circ L_2'.$

Let ${\cal F}_2'= L_2({\cal F}_2)\cup \{y_1,y_2\}.$
Let ${\cal G}_2={\cal G}(\ep_3/2, {\cal F}_2'),$
${\cal Q}_2={\cal P}(\ep_3/2, {\cal F}_2')$ and
$d_2=\dt(\ep_3/2, {\cal F}_2')$ be as in \ref{IIIL3} (for $B$)
associated with $\ep_3/2$ and ${\cal F}_2'.$
We may assume that ${\cal F}_2'\supset {\cal F}_2',$
${\cal G}_2 \supset [L_2]({\cal P}_2)$ and
$d_2<{\rm min}\{\ep_3/2, \dt_2/2\}.$

It follows from \ref{IIIL2} that there is a \morp\,
$\Psi_2': B\to A$ such that
$$
\|\Psi_2'(cd)-\Psi_2'(c)\Psi_2'(d)\|<d_2/2
$$
for all $c,d\in {\cal G}_2$ and
$$
[\Psi_2']|_{{\cal Q}_2}=\alpha^{-1}|_{{\cal Q}_2}.
$$
By \ref{IIIL3} there is a unitary $u_2\in A$ such
that
$$
\Ad u_2\circ \Psi_2'\circ L_2\approx_{\ep_3/2} \id_A
\,\,\,\,\,\,{\rm on}\,\,\, {\cal F}_2.
$$
Set $\Psi_2=\Ad u_2\circ \Psi_2.$

Continuing in this fashion, we construct a
sequence of \morp s $L_n: A\to B$ and $\Psi_n: B\to A$ such that
the following diagram
$$
\begin{array}{ccccccc}
\noalign{\vskip-9pt}
A  & {\stackrel{\id_A}{\longrightarrow}}  & A
&{\stackrel{\id_A}{\longrightarrow}}
& A& {\stackrel{\id_A}{\longrightarrow}} &\cdots A\\
\downarrow_{L_1}& \nearrow_{\Psi_1} &\downarrow_{L_2}& \nearrow_{\Psi_2}&
\downarrow_{L_3}& \nearrow_{\Psi_3} &\\
B&{\stackrel{\id_A}{\longrightarrow}}& B&
{\stackrel{\id_A}{\longrightarrow}}& B&
{\stackrel{\id_A}{\longrightarrow}}&\cdots B\\
\end{array}
$$
is approximately intertwining.
It follows from an argument of Elliott (see for example 2.1, 2.2 and 2.3 in
\cite{Ell1}
and also 3.1 in \cite{Ln1} for the case that the maps are not \hm s)
that there are isomorphisms $h: A\to B$ and
$h^{-1}: B\to A$ (each determined by $\{L_n\}$ and $\{\Psi_n\}$).
\enddemo

\proclaim{Theorem}\label{IIIT2}
Let $A$ and $B$ be two unital \CA s which are inductive
limits of \CA s in ${\cal BD}$ with $ {\rm TR}(A)= {\rm TR}(B)=0.$
Suppose that there is an order isomorphism $$\alpha:
(K_0(A), K_0(A)_+, [1_A], K_1(A))\to
(K_0(B), K_0(B)_+, [1_B], K_1(B));$$ 
then there is an isomorphism $h: A\to B$ such
that $h_*=\alpha.$
\endproclaim 

\demo{Proof}
Since each \CA\, is in ${\cal BD},$ both $A$ and $B$ satisfy
the UCT. So the theorem follows from \ref{IIIT1}.
\enddemo

Recently, it was shown by Q. Lin and N.C. Phillips
that the following structure theorem about
smooth minimal dynamical systems holds:

\proclaimtitle{Q. Lin and N.C. Phillips \cite{LP}}
\proclaim{Theorem}
\label{IIILP3}
Let $M$ be a compact manifold and $\dt: M\to M$ be a minimal
diffeomorphism. Then $A=C^*({\Bbb Z}, M,\dt),$ the simple
crossed  product arising from the smooth minimal dynamical system{\rm ,} is
in ${\cal LBD}.$ In fact{\rm ,} $A$ is a direct limit
of subhomogeneous \CA s.
Furthermore $A$ has stable rank one and $K_0(A)$ is
weakly unperforated.
\endproclaim 

Thus we have the following:

\proclaim{Theorem}\label{IIIT5}
Let $M_1$ and $M_2$ be  compact manifolds, $\dt_i: M_i\to M_i$
be a minimal
diffeomorphism {\rm (}$i=1,2${\rm )} and let $A_i=C^*({\Bbb Z}, M_i,h_i).$
Suppose that $ {\rm TR}(A)=0.$
Then $A_1\cong A_2$ if and only if
$$
(K_0(A_1), K_0(A_1), [1_{A_1}], K_1(A_1))\cong
(K_0(A_2), K_0(A_2)_+, [1_{A_2}], K_1(A_2)).
$$
\endproclaim 

An important case for dynamical systems is the unique ergodic case,
where systems admit unique invariant measures.
The resulting simple crossed products
admit unique normalized traces.

\proclaim{{C}orollary}\label{CROSS3}
Let $M_1$ and $M_2$ be compact manifolds,
$h_i: M_i\to M_i$ be a minimal diffeomorphism {\rm (}$i=1,2${\rm ).}
Let $A_i=C^*({\Bbb Z}, M_i, h_i).$
Suppose that $(M_i, h_i)$ has a unique invariant measure {\rm (}$i=1,2${\rm )}
and the range of $K_0(A_i)$ under the trace is dense in
${\Bbb R}.$
If
$$
(K_0(A_1), K_0(A_1), [1_{A_1}], K_1(A_1))\cong
(K_0(A_2), K_0(A_2)_+, [1_{A_2}], K_1(A_2)),
$$
then
$$
A_1\cong A_2.
$$
\endproclaim

\demo{Proof}
By \cite{LP}, the  $A_i$ are simple \CA s that have stable rank one
and have weakly unperforated $K_0(A_i).$
It follows from \cite{Ph}
that both $A_i$ have real rank zero.
Thus the corollary follows from \ref{UNI} and \ref{IIIT5}.
\enddemo

There is a class of exciting unital simple
\CA s   called (higher) irrational noncommutative tori.
The irrational rotation \CA\, $A_{\theta}$ is generated by
two unitaries $u$ and $v$ with relation $uv=e^{i\pi\theta}vu,$
where $\theta$ is irrational, and is a very interesting
unital simple \CA.
It was shown that it is a direct limit of circle algebras
with real rank zero (see \cite{EE}). Therefore
$ {\rm TR}(A_{\theta})=0.$ Also $A_{\theta}$ can  be realized as
the crossed product $C({\Bbb T})\times_{\theta} {\Bbb Z}$
resulting from the minimal and unique ergodic
dynamical system by a diffeomorphism $\alpha:
{\Bbb T}\to {\Bbb T}$ which
maps $z$ to $ze^{i\pi\theta},$ where $\theta$ is irrational.
One class of $k+1$-dimensional noncommutative tori
called unital simple crossed products
$C({\Bbb T}^k)\times_{\theta}{\Bbb Z}$
result from  a minimal and unique ergodic
dynamical system by a diffeomorphism $\alpha:
{\Bbb T}^k\to {\Bbb T}^k$ which
maps $(z_1,\ldots ,z_k)$ to $(z_1e^{i\pi\theta_1},\ldots ,z_ke^{i\pi\theta_k}),$
where $\theta_i$ are irrational.

We have the following theorem:

\proclaim{Theorem}\label{IVT2}
Let $A=C({\Bbb T}^k)\times_{\theta} {\Bbb Z}$
and $B=C({\Bbb T}^k)\times_{\alpha} {\Bbb Z}$ be  two
irrational noncommutative $k+1$\/{\rm -}\/tori.
Then $A\cong B$ if and only if
$[\alpha]=[\theta]$ on $K_i(C({\Bbb T}^k))$ {\rm (}$i=0,1$ {\rm )}.
Furthermore{\rm ,} $A$ is an inductive limit
of circle algebras.
\endproclaim

\demo{Proof}
By \cite{PV}, $K_i(A)=K_i(B)$
($i=0,1$).
It follows from \ref{CROSS3}
that $A\cong B.$ They are also isomorphic to a unital simple \CA\,
in ${\cal C}.$ Since $K_i(A)$ are torsion free ($i=0,1$),
they are isomorphic to an inductive limit of circle algebras.

\numbereddemo{{R}emark}
{\rm One should note that the condition that $A_m$ are
subhomogeneous in \ref{IIIL2} is not necessary.
It suffices to have the following:
any finitely generated subgroup $G\subset \rho_{A_n}(K_0(A))$ can
be embedded into ${\Bbb Z}^k$ for some integer $k>0$ as
an ordered subgroup. It remains open whether this condition
is automatically satisfied by residually finite-dimensional
algebras, or any simple nuclear \CA s which can be written as
inductive limits of RFD \CA s satisfying this condition.
}
\enddemo
 
\AuthorRefNames [EHS]

\end{document}